\documentclass{sig-alternate}

\usepackage{todonotes}

\usepackage{overpic}
\usepackage{amsmath,amstext,amsbsy,amssymb,subfig,mathdots}
\usepackage{algorithm,algpseudocode}
\usepackage{booktabs}
\usepackage{mathrsfs}
\usepackage{default,tikz}
\usetikzlibrary{positioning}
\usetikzlibrary{shapes,arrows}

\newenvironment{remark}{{\bf Remark}}

\newtheorem{example}{Example}

\def\L{{\cal L}}
\def\O{{\cal O}}
\def\M{{\cal M}}
\def\S{{\cal S}}

\def\P{{\cal P}}

\def\dkf{{\mathrm d}}

\definecolor{DarkGreen}{rgb}{0,.55,0}

\definecolor{Yellow}{rgb}{.55,.55,0}

\let\oldremark\remark\renewcommand{\remark}{\oldremark\normalfont}
\let\oldexample\example\renewcommand{\example}{\oldexample\normalfont}

\newcount\examplenumber
\def\Example#1{\advance \examplenumber by 1
\subsection{Example \the\examplenumber\ (#1)}}

\def\Ki_#1^#2{{\mathscr K}_{#1}^{#2}}

\def\sopmatrix#1{\begin{pmatrix} #1 \end{pmatrix}}

\def\figref#1{Figure~\ref{fig:#1}}

\def\addtab#1={#1\;&=}

\def\meeq#1{\def\ccr{\\\addtab}
 \begin{align*}
 \addtab#1
 \end{align*}
  }

\def\vc#1{\mbox{\boldmath$#1$\unboldmath}}

\def\Figure#1#2\par{
\begin{figure}[tb]
\begin{center}{
\includegraphics{Figures/#1}}
\end{center}
\caption{#2}\label{fig:#1} 
\end{figure}
}
\def\Figurew#1#2#3\par{
\begin{figure}[tb]
\begin{center}{
\includegraphics[width=#2]{Figures/#1}}
\end{center}
\caption{#3}\label{fig:#1} 
\end{figure}
}

\def\Figuretwow#1#2#3#4\par{
\begin{figure}[tb]
\begin{center}{
\includegraphics[width=#3]{Figures/#1}\includegraphics[width=#3]{Figures/#2}}
\end{center}
\caption{#4}\label{fig:#1} 
\end{figure}
}
\def\Figuretwo#1#2#3\par{
	\Figuretwow{#1}{#2}{0.48 \hsize}
		#3\par	
}

\def\pr(#1){\left({#1}\right)}
\def\br[#1]{\left[{#1}\right]}
\def\set#1{\left\{{#1}\right\}}
\def\ip<#1>{\left\langle{#1}\right\rangle}
\def\iip<#1>{\left\langle\!\langle{#1}\right\rangle\!\rangle}

\def\fpr(#1){\!\pr({#1})}

\def\function#1#2{\expandafter\def\csname #1\endcsname(##1){#2\fpr({##1})}
				\expandafter\def\csname #1p\endcsname(##1){#2'\fpr({##1})}
				\expandafter\def\csname #1pp\endcsname(##1){#2''\fpr({##1})}
				\expandafter\def\csname #1pn\endcsname##1(##2){#2^{\pr({##1})}\fpr({##2})}				}
\def\defoperator#1#2{\expandafter\def\csname #1\endcsname[##1]{#2\!\br[{##1}]}}

\def\ffunction#1{\function{#1}{#1}}

\def\O(#1){{\cal O}\!\left(#1\right)}

\def\Oo(#1){{\rm o}\!\left({#1}\right)}

\def\fO(#1){\Oh\fpr({#1})}

\def\dkfn^#1#2{\,{\rm d}^#1 #2}
\def\dkf#1{\,{\rm d}#1}

\def\H_#1^#2(#3){H_#1^{(#2)}\fpr({#3})}
\def\J_#1(#2){{\rm J}_#1\!\pr(#2)}

\def\dz{\dkf{z}}

\def\ddxn^#1{{\dkf{}^#1 \over \dkfn^#1{x}}}

\def\seq_#1^#2#3{\set{#3_{#1},\ldots,#3_{#2}}}

\ffunction{f}
\ffunction{g}

\def\mapengine#1,#2.{\mapfunction{#1}\ifx\void#2\else\mapengine #2.\fi }

\def\map[#1]{\mapengine #1,\void.}

\def\mapenginesep_#1#2,#3.{\mapfunction{#2}\ifx\void#3\else#1\mapengine #3.\fi }

\def\mapsep_#1[#2]{\mapenginesep_{#1}#2,\void.}

\def\vcbr[#1]{\pr({#1})}

\def\bvect[#1,#2]{
{
\def\dots{\cdots}
\def\mapfunction##1{\ | \  ##1}
	\sopmatrix{
		 \,#1\map[#2]\,
	}
}
}

\def\vect[#1]{
{\def\dots{\ldots}
	\vcbr[{#1}]
}}

\def\vectt[#1]{
{\def\dots{\ldots}
	\vect[{#1}]^{\top}
}}

\def\Vectt[#1]{
{
\def\mapfunction##1{##1 \cr} 
\def\dots{\vdots}
	\sopmatrix{
		\map[#1]
	}
}}

\def\qqqquad{\qquad\qquad}

\def\qfor{\quad\hbox{for}\quad}

\def\simlimit_#1{\,\,\,\sim \!\!\!\!\!\!\!\!\!{ \atop \scriptscriptstyle #1 }}

\def\XXint#1#2#3{{\setbox0=\hbox{$#1{#2#3}{\int}$}
     \vcenter{\hbox{$#2#3$}}\kern-.5\wd0}}

\def\rad^#1#2{\,\,{}^#1\!\!\!\!\sqrt{#2}\,}

\def\Figuretwofixed#1#2#3\par{
\Figuretwow{#1}{#2}{0.48 \hsize}{#3}\par
}

\title{A practical framework for infinite-dimensional linear algebra}

\author{Sheehan Olver\thanks{School of Mathematics and Statistics, The University of Sydney, Sydney, Australia. (Sheehan.Olver@sydney.edu.au)} 
\and 
Alex Townsend\thanks{Department of Mathematics, Massachusetts Institute of Technology, 77 Massachusetts Avenue
Cambridge, MA 02139-4307. (ajt@mit.edu)} 
}

\begin{document}
\maketitle

\begin{abstract}
We describe a framework for solving a broad class of infinite-dimensional linear equations, consisting of almost banded operators, which can be used to resepresent linear ordinary differential equations with general boundary conditions.   The framework contains a data structure on which row operations can be performed, allowing for the solution of  linear equations by the adaptive QR approach.  The algorithm achieves $\O(n^{\rm opt})$ complexity, where $n^{\rm opt}$ is the number of degrees of freedom required to achieve a desired accuracy, which is determined adaptively.  In addition, special tensor product equations, such as  partial differential equations on rectangles, can be solved by truncating the operator in the $y$-direction with $n_y$ degrees of freedom and using a generalized Schur decomposition to upper triangularize, before applying the adaptive QR approach to the $x$-direction, requiring $\O(n_y^2 n_x^{\rm opt})$ operations.  The framework is implemented in the {\sc ApproxFun} package written in the {\sc Julia} programming language, which achieves highly competitive computational costs by exploiting unique features of {\sc Julia}.
\end{abstract}


%\begin{AMS}
%33A65, 35C11, 65N35
%\end{AMS}

\section{Introduction}

	Linear equations play a fundamental role in scientific computing, with the classical examples including the numerical solution of  boundary value ordinary differential equations, elliptic partial differential equations and singular integral equations.    Practically all numerical methods for solving linear differential equations --- e.g., finite difference, finite element, collocation and Galerkin methods --- can be described as {\it discretize-then-solve}. That is, the underlying infinite-dimensional operator is first approximated by a finite-dimensional matrix, before the resulting linear system is solved by a standard linear algebra method that is either direct, such as Gaussian elimination, or iterative, such as conjugate gradient.   
	
		In contrast, we advocate an entirely different approach that solves the equation as an infinite-dimensional problem and never discretizes the operator itself.  To accomplish this task, we represent the (infinite-dimensional) operator by a suitable data structure that supports  row manipulations directly on the representation of the operator, using lazy evaluation to automatically extend the data in the representation as needed.  The row operations can be used to partially upper-triangularize the operator, and, for a large class of problems (in particular non-singular ODEs), we can at some point  perturb the right-hand side by a small amount so that the still infinite-dimensional problem can be solved exactly via a (finite-dimensional) back substitution step.
		
			The mathematical ground work for this approach is (F. W. J.) Olver's algorithm~\cite{OlversAlgorithm}, which considers the solution of inhomogeneous three-term recurrence relationship.
%%
%	\meeq{
%		b_0 u_0 + b_1 u_1 + b_2 u_2 + \cdots = f_0, \ccr
%		\gamma_r u_{r-1} + \alpha_r u_r + \beta_r u_{r+1} =f_r \qfor r=1,2,3,\ldots.
%		}
%%
This is equivalent to solving an infinite-dimensional linear system involving a rank-1 perturbation of a tridiagonal operator:
	$$\sopmatrix{b_0 & b_1 & b_2 & b_3 & \cdots \\
				\gamma_1 & \alpha_1 & \beta_1 \\
				& \gamma_2 & \alpha_2 & \beta_2 \\
				&& \gamma_3 & \alpha_3 & \ddots \\
				&&&\ddots & \ddots } \begin{pmatrix}u_0\\[2pt] u_1 \\[2pt] u_2\\[2pt] u_3\\[2pt] \vdots\end{pmatrix} = \begin{pmatrix} f_0\\[2pt] f_1\\[2pt]f_2\\[2pt]f_3\\[2pt]\vdots\end{pmatrix}.$$
The key observation is that the infinite-dimensional linear system can be solved by Gaussian elimination without pivoting and that convergence to the {\it minimal solution} of the system --- roughly, the solution (provided it exists) with the fastest decaying entries ---   can be inferred as part of the algorithm.  
Back substitution then proceeds by perturbing the right-hand side, as opposed to changing the infinite-dimensional operator.  In the functional analysis setting, where the operator is assumed to be invertible between two spaces, the minimal solution is the unique solution to the linear equation.  
This approach was extended by Lozier to more general banded operators~\cite{Lozier_80_01}. However, Gaussian elimination without pivoting is prone to numerical instability,  and with this in mind the authors derived an adaptive QR approach~\cite{Olver_13_01}, using Givens rotations for the solution of linear ordinary differential equations (ODEs).  The complexity of Olver's algorithm and the adaptive QR approach is $\O(n^{\rm opt})$, where $n^{\rm opt}$ is the number of coefficients calculated, as determined automatically by the convergence criteria.  (Throughout, an integer with a superscript ``opt'' is a number that is determined adaptively as dictated by the particular problem.)

Similar in spirit to the current work is Hansen~\cite{Hansen_08_01}, which investigated the infinite-dimensional QR algorithm for spectral problems, though focusing on theoretical rather than practical matters.   Our operator algebraic framework is heavily influenced by the {\tt chebop} system~\cite{Driscoll_08_01}, which is part of {\sc Chebfun}~\cite{Chebfun}, and provides an infinite-dimensional feel to the user, though the underlying collocation method is the traditional approach of discretize-then-solve.   Finally, the second author and Trefethen investigated  continuous analogues of matrix algorithms~\cite{Townsend_14_01}, where the emphasis is on representing smooth bivariate functions rather than operators.

	In this work, we exploit the applicability of the adaptive QR method for a general class of linear operators: banded operators except for possibly a finite number of dense rows.  However, to make this competitive and useful for general problems requires the following components:
\begin{enumerate}
	\item Abstract data types that can be overriden to represent arbitrary (typically unbounded) banded operators and dense functionals, as well as data structures to allow their algebraic manipulation.
	\item {\it Fast} linear algebra on infinite-dimensional operators. This requires a carefully managed data structure that can encapsulate the full infinite-dimensional operator at each stage of the linear algebra routines.
\end{enumerate}
Unfortunately, the specialized data structures that we develop also require 
very specific implementations of linear algebra routines,  prohibiting the 
traditional approach of reducing the problem to finite-dimensional linear algebra solvable by {\sc LAPack}.   In~\cite{Olver_13_01}, the C++ language was used to partially  implement the framework for some simple examples; however, adding new operators required a complete recompilation, which is prohibitively time consuming for practical use.  As an alternative, the {\sc Julia} programming 
language~\cite{Julia} provides a natural 
environment for implementing both data structures and linear 
algebra algorithms.  Furthermore, the support for multiple dispatch allows for the 
easy construction of data structures, and linear algebra can be performed remarkably efficiently
due to on-the-fly compilation. The {\sc ApproxFun} package \cite{ApproxFun} implements the proposed framework in {\sc Julia}.

%In~\secref{DataStructure} we introduce the general framework for infinite-dimensional operators.  We give a simple toy example of differential operators on Taylor coefficients; for the construction of real world differential operators that fit into the framework,  we refer the reader  to the operators developed in the ultraspherical spectral method~\cite{Olver_13_01}.  In \secref{LinAlg} we describe a wrapper data structure that supports linear algebra on infinite-dimensional operators, leading immediately to the solution of linear systems of equations.  In \secref{PDEs}, we see that this can be combined with a discretization in one  dimension to solve tensor product problems such as partial differential equations a l\'a~\cite{Townsend_14_03}.    By using the infinite-dimensional approach for the remaining dimension, we achieve faster computational times, and indeed better complexity.

\section{Data structures for  operators}\label{sec:DataStructure}
We represent operators  by one of the following fundamental abstract types:
\begin{enumerate}
\item {\tt Functional}: A data structure representing an operator of size $1 \times \infty$. 
\item {\tt BandedOperator}: A data structure representing an operator of size $\infty \times \infty$ that has a finite bandwidth with the bands ranging from $a\!:\!b$, where $a\leq 0 \leq b$. That is, the $k$th row only has (possibly) nonzero entries in columns $a+k,\ldots,b+k$. 
\end{enumerate}
A new functional, say {\tt NewFunctional}, which is a subtype of {\tt Functional}, must override a routine called {\tt getindex}, 
$$\hbox{\tt getindex(F::NewFunctional,cr::Range)},$$
which returns a vector of the entries in the columns specified by {\tt cr}.  Similarly, each subtype of {\tt BandedOperator} overrides a routine called {\tt addentries!}, which adds entries to specified rows of a (finite) banded array\footnote{The use of the ``!''  suffix  is a {\sc Julia} convention that signifies a method that modifies one of its inputs, in this case, a banded array.  We add, as opposed to overwrite, to make addition of operators more efficient.  The details of the banded array data structure are immaterial, but an $n \times n$ banded array with bands ranging from $a\!:\!b$ can be represented by an $n \times (b-a+1)$ array.}, and overrides a routine called {\tt bandinds} that returns the band range of the operator represented as a tuple $(a,b)$.  

As an example, consider representing functionals and operators that act on vectors of Taylor series coefficients, i.e., vectors of the form $\vectt[u_0,u_1,u_2,\dots]$ that correspond to the series $\sum_{k=0}^\infty u_k z^k$.  Evaluation at a point $z$ is the functional ${\cal B}_z \triangleq [1,z,z^2,\dots]$ and thus we can create a subtype {\tt TaylorEvaluation}, with single field {\tt z}, that implements 
	$$\hbox{\tt getindex(B::TaylorEvaluation,cr)=B.z.\textasciicircum(cr-1)}.$$
We can also implement a {\tt TaylorDerivative} operator to represent the banded operator, defined by $\vc e_k^\top {\cal D} \vc e_{k+1} = k$ and zero otherwise, 
%%
%	$${\cal D}\triangleq\sopmatrix{0 & 1 \cr &&2 \cr &&& 3 \cr &&&&\ddots}$$
%%
with band range
%\footnote{The implementation currently requires a band range $a\!:\!b$ where $a\leq0$ and $b\geq0$. Otherwise, the band range could be taken as $1\!:\!1$.}  
$0\!:\!1$.  Finally, we can represent multiplication by a polynomial of finite degree, says $a(z)=\sum_{k=0}^{m-1} a_k z^k$ with the Toeplitz operator defined by $\vc e_k^\top {\cal T}[a]\vc e_j = a_{k-j}$ for $0\leq k-j \leq m-1$, and zero otherwise.  
%	$${\cal T}[a]\triangleq\sopmatrix{a_0 \cr a_1 & a_0 \cr a_2 & a_1 & a_0 \cr \vdots & \ddots & \ddots & \ddots   \cr a_{m-1}  &\cdots & a_2 & a_1  & a_0 \cr  & \ddots & \ddots & \ddots & \ddots   & \ddots }.$$  
	This is encoded in a {\tt TaylorMultiplication} operator, which has a single field containing the coefficients of $a$ as a vector of length $m$ and a band range $(1-m)\!:\!0$.

	While our operators always act on infinite-dimensional vectors, the entries of those vectors can represent coefficients in many different bases. For example, in~\cite{Olver_13_01} vectors represent expansion coefficients in the Chebyshev or ultraspherical basis.  To ensure that the domain and range of two operators are consistent when, for instance, adding them together, each operator must know the basis of its domain and range. This also allows us to automatically convert between bases to ensure that any operation can be performed in a consistent manner.  Therefore, we have a {\tt  FunctionSpace} abstract type so that operators can override {\tt domainspace} and {\tt rangespace} routines that return specific domain and range spaces.  When it exists, a banded conversion operator is implemented to convert between two spaces.  
	
	The final components are structures that allow {\it functional} and {\it operator} algebra.  This consists of a {\tt PlusOperator}, which contains a list of {\tt BandedOperator}s that have the same range and domain spaces.  The command ``$+$'' is then overridden for {\tt BandedOperator}s, with an additional step of promoting the domain and range space whenever a banded conversion operator is available.  Similarly, a {\tt TimesOperator} is constructed to represent multiplication of operators, and ``$*$'' is similarly overloaded to promote spaces to ensure compatibility.  Note that, if {\tt A} and {\tt B} have band range $a\!:\!b$ and $c\!:\!d$, respectively, then the band range of {\tt A * B} is $(a+c)\!:\!(b+d)$.  To determine the entries of {\tt A * B} up to row $k$, we determine {\tt A} up to row $k$, {\tt B} up to row $k+b$ and multiply as appropriate.  
 Similarly, a {\tt PlusFunctional} is implement to represent addition of {\tt Functional}s and a {\tt TimesFunctional} to represent a {\tt Functional} times a {\tt BandedOperator}.

 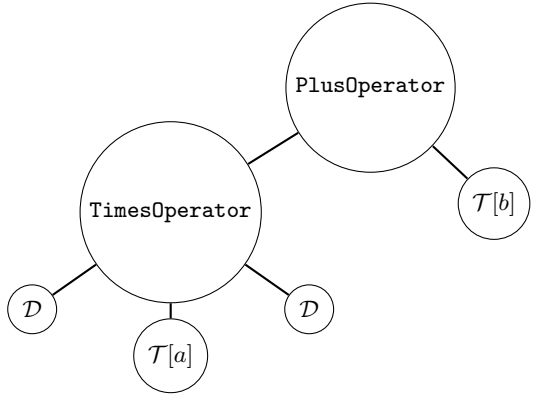
\begin{figure}
\centering
\tikzstyle{block} = [circle, draw, fill=white, minimum width=2em]
\tikzstyle{line} = [draw, thick, color=black]
\hspace*{18pt}\begin{tikzpicture}[scale=2, node distance = 1cm, auto]
    % Place nodes
    \node [block] (top) {{\tt PlusOperator}};
    \node [block, below left = 0cm and 1cm of top] (mult) {{\tt TimesOperator}};
    \node [block, below right = .4cm and .5cm of top] (b) {${\cal T}[b]$};
    \node [block, below left = .2cm and .75cm of mult] (yseed) {${\cal D}$};
    \node [block, below = .2cm and .5cm of mult] (seed6) {${\cal T}[a]$};    
    \node [block, below right = .2cm and .75cm of mult] (seed4) {${\cal D}$};
    % Draw edges
    \path [line] (top) -- (mult);
    \path [line] (top) -- (b);
    \path [line] (mult) -- (yseed);
    \path [line] (mult) -- (seed4);
    \path [line] (mult) -- (seed6);    
\end{tikzpicture}
\caption{ 	General banded operators can be built up from elementary operators using {\tt PlusOperator} and {\tt TimesOperator}.   This tree represents the operator ${\dkf \over \dz} a(z) {\dkf \over \dz} + b(z)$ being built up from differentiation operator $\mathcal{D}$ and multiplication operators $\mathcal{T}$.
}
\label{fig:OperatorTree}
\end{figure}

 	As operators are manipulated algebraically, a tree structure is automatically constructed.  Returning to the Taylor series example, we can represent the differential operator ${\dkf \over \dz} a(z) {\dkf \over \dz} + b(z)$ as ${\cal D} {\cal T}[a] {\cal D} + {\cal T}[b]$, which has the tree structure as depicted in \figref{OperatorTree}.  
	We expect the tree to be small in most cases, providing  a {\tt SavedOperator} type that wraps a banded operator to save its entries as they are computed.
	We finally mention that there is an interlace operator that takes two or more operators and alternates their entries.  This facilitates a natural extension to (small) systems of differential equations, as well as problems posed on multiple domains, where the continuity conditions are represented as functionals.  This also allows us to work with doubly-infinite operators (e.g., operators acting on Fourier series), via interlacing the non-negative and negative entries.

\section{Infinite-dimensional linear algebra}\label{sec:LinAlg}

	Using the above structure, we can represent quite general operators with boundary conditions by a list of $K$ {\tt Functional}s and a single {\tt BandedOperator}. For example,  
	$${\dkf \over \dz} a(z) {\dkf u \over \dz} + b(z) u = f(z),u(1)=c_1, u'(1)=c_2$$
	 represented with Taylor series becomes
	$$\Vectt[{\cal B}_1, {\cal B}_1{\cal D} , {\cal D} {\cal T}{[a]} {\cal D} + {\cal T}{[b]}]  u =  \Vectt[c_1,c_2,f].$$
Here, the first entry is a {\tt TaylorEvaluation}, the second entry is a {\tt TimesFunctional} with leaf nodes {\tt TaylorEvaluation} and {\tt TaylorDerivative}, and the third entry is a {\tt PlusOperator}, which is the root node of the tree structure depicted in \figref{OperatorTree}.  
This is an {\it almost-banded operator}, in the sense that it is a banded operator (with band range $(a-K)\!:\!(b-K)$) except for the first $K$ dense rows.

	We now wrap this operator by a {\tt MutableAlmostBandedOperator}, a mutable data structure that allows for row operations, i.e., the addition of one row to the other.    In this case, the represented operator is an almost-banded operator where the first $(K+n)$ rows are dense and $n$ is an integer that is selected adaptively. The remaining rows have band range $(a - K) \!:\! (b-a)$. More specifically, its fields are as follows: 
\begin{enumerate}
	\item {\tt F}: A  list of $K$ {\tt Functional}s.
	\item {\tt B}: A single {\tt BandedOperator} of band range $a\!:\!b$.
	\item {\tt bcdata}: A $K \times  (b - a + K)$ array containing the mutable entries of the boundary rows, where only the entries on and below the $b-a$ super diagonal are used.
	\item {\tt bcfilldata}: A $K \times K$ array which dictates how the first $K$ rows are filled in.
	\item {\tt data}: An $n \times n$ banded array with band range $a\!:\!(b-a+K)$ containing the mutable entries of the banded operator.
	\item {\tt filldata}: A  $n \times K$ array which dictates how the $(K+1)$th through $(K+n)$th rows are filled in.
\end{enumerate} 
If $A$ is a {\tt MutableAlmostBandedOperator}, then it represents the $\infty \times \infty$ operator with the $k,j$th entry given by
	$$A[k,j]  = \begin{cases}
		{\tt bcdata}[k,j], \cr \qqqquad \qfor 1 \leq k \leq K, 1\leq j < M+k, \cr
		\sum_{i=1}^K {\tt bcfilldata}[k,i] {\tt F}[i][j], \cr
				\qqqquad	 \qfor 1 \leq k \leq K,   M+k  \leq j, \cr
		{\tt data}[k-K,j], \cr
		 \qqqquad	 \qfor K < k \leq n+K, 1 \leq j < M+k, \cr
		\sum_{i=1}^K {\tt filldata}[k-K,i] {\tt F}[i][j], \cr
		\qqqquad	 \qfor K < k \leq n+K,  M+k \leq j, \cr
		{\tt B}[k-K,j], \cr
		\qqqquad	 \qfor  n+K < k, \cr
	\end{cases}$$
where $M=b-a + K$.
	
	We perform row operations that column-by-column introduce zeros below the diagonal. Therefore, when acting on rows $k_1$ and $k_2$ (with $k_1 < k_2$) the first $k_1-1$ entries in both rows are already zero.   Under such circumstances, $K$ remains fixed and the complexity of introducing zeros in the first $n$ columns is $\O((K-a)^2 n)$, with an $\O(n)$ growth of data storage.    We refer the reader to \cite{Olver_13_01} for a precise description on how this is achieved using Givens rotations. 
	
%	If we act on rows $k_1$ and $k_2$ by the matrix $\sopmatrix{a & b \cr c & d}$, then we update the data according to:
%%
%\sotodoinline{Describe procedure precisely}
%%
%	
	For solving linear equations, we also apply row operations to the right-hand side.  If the initial right-hand side has a finite number of nonzero entries, then this can be done adaptively.   Suppose that, after introducing zeros below the diagonal in the first $n$ columns, the right-hand  side happens to have has zeros apart from its first $n$ entries.  Then it lies in the span of the upper triangular component of $A$, and we can proceed with back-substitution, in $\O(n)$ operations, see \cite{Olver_13_01}.  
If the entries past the $n$th entry are small, we can truncate them to produce a new right-hand side, close to the original right-hand side, so that the resulting equation is solvable by back substitution.  We emphasize this is a distinct process from changing the operator: we can control the effect of changing the right-hand side, and avoid issues of causing an invertible operator to become non-invertible.   We expect well-posed problems to converge via this methodology, see \cite{Olver_13_01} for a proof in the case of non-singular ordinary differential equations.

Other linear algebra routines are also applicable in infinite dimensions using row manipulations on the same data structure.  For example,  the null space of a banded operator can be calculated by applying Givens rotations to the transpose of the operator, see the {\tt null} command in {\tt ApproxFun} \cite{ApproxFun}.    Still under investigation is whether the infinite-dimensional QL algorithm\footnote{In finite dimensions, the QL and QR algorithms are equivalent.  In infinite dimensions this is no longer the case, and it is  likely that only an infinite-dimensional QL algorithm can support shifts to induce faster convergence.  The infinite-dimensional QL algorithm does not appear to have been investigated, unlike the the infinite-dimensional QR algorithm \cite{Hansen_08_01} and the infinite-dimensional Toda flow \cite{Deift_85_01}. }   with Wilkinson shifts can be implemented for calculating spectrum of operators.  

\section{Infinite-dimensional tensor equations}\label{sec:PDEs}

	We finally consider general tensor equations with two terms --- that is, {\it splitting rank} 2 in the terminology of \cite{Townsend_14_03} ---  with boundary conditions.  We represent such equations as acting on an unknown $\infty \times \infty$ matrix $X$ that satisfies:
\def\N{{\cal N}}
	\meeq{
		\L X \M^\top + \N X \S^\top = F,\qquad X{\cal B}_y^\top = \vc g_y,\qquad {\cal B}_x X = \vc g_x^\top,}
%	\meeq{
%		({\cal B}_x \otimes I) u = \vc g(y)\qquad (I \otimes {\cal B}_y) u = \vc h(x)\qqand(\L \otimes \M + \N \otimes \S) u = f}
%
where $F$ is an $\infty \times \infty$ matrix corresponding to a forcing term, $\vc g_x$ and $\vc g_y$ are $\infty \times K_x$ and $\infty \times K_y$ matrices corresponding to boundary condtions, ${\cal B}_x$ and ${\cal B}_y$ are vectors of $K_x$ and $K_y$ functionals, and $\L,\M,\N$ and $\S$ are banded operators.  Many standard linear PDEs on rectangles with boundary conditions  --- e.g., Helmholtz equation, Poisson equation, linear KdV and the semi-classical Schr\"odinger equation with a time-independent potential ---  can be written in this form using the ultraspherical spectral method to obtain banded differential operators~\cite{Townsend_14_03}.   
%We assume that the equation has a unique solution under suitable conditions, and do not address the convergence of the proposed algorithm to such solution\attodo{Why do we need to say this?}.  

% If we represent $u$ and $f$ as $\infty \times \infty$ matrices $X$ and $F$, this equation takes the form
%

%
	The  equation we wish to solve is an infinite-dimensional analogue of a {\it generalized Sylvester equation}~\cite{Gardiner_92_01}.  We adapt the approach of~\cite{Townsend_14_03}, which solved this equation by  discretizing and upper triangularizing in both dimensions, to now  only discretizing in one dimension.  Define the projection operator ${\cal P}_n:\mathbb{C}^\infty\rightarrow\mathbb{C}^n$ as 
	$${\cal P}_n\vectt[u_0,u_1,\dots]=\vectt[u_0,u_1,\dots,u_{n-1}],$$
	and  consider $X_n$, the $\infty \times n(\equiv n_y)$ solution to the semi-discretized equation:
	\meeq{
	      \L X_n M_n^\top + \N X_n S_n^\top = F_n, \ccr
	       {\cal B}_x X_n = \vc g_{xn}^\top, \qquad X_nB_{n}^\top = \vc g_y,}
where $M_n = {\cal P}_n \M {\cal P}_n^\top$,  $S_n = {\cal P}_n \S {\cal P}_n^\top$, $B_{n} = {\cal B}_y {\cal P}_n^\top$, $\vc g_{xn} = {\cal P}_n \vc g_x$ and $F_n = F {\cal P}_n^\top$.  

Assume without loss of generality that  $B_n = \bvect[I_{K_y},B_n^{(2)}]$, i.e.,  the principle $K_y \times K_y$ block of $B_{n}$ is the identity matrix, see~\cite{Townsend_14_03} for the procedure to ensure that this is true.   We  incorporate the  discretized boundary conditions into the generalized Sylvester equation by removing the dependence on the first $K_y$ columns of $X_n$, i.e., introducing zeros in the first $K_y$ rows of $M_n^\top$ and $S_n^\top$ via
		\meeq{
			\L X_n (M_n^\top - B_n^\top M_{n,K_y}^\top) + \N X_n (S_n^\top - B_n^\top  S_{n,K_y}^\top)  = \cr
			F_n - \L \vc g_y - \N \vc g_y,
			}
where $M_{n,K_y}$ and $S_{n,K_y}$ are the $n \times K_y$ principle subblocks of $M_n$ and $S_n$, respectively.  

Partitioning $X_n = \bvect[X_n^{(1)}, X_n^{(2)}]$ so that $X_n^{(1)}$ is $\infty \times K_y$ and $X_n^{(2)}$ is $\infty \times (n-K_y)$, we see that $X_n^{(2)}$ satisfies
	$$\L X_n^{(2)} \tilde M_n^\top + \N X_n^{(2)} \tilde S_n^\top = \tilde F_n$$
for suitable matrices $\tilde M_n$, $\tilde S_n$ and $\tilde F_n$.   We now modify the Bartels--Stewart algorithm~\cite{Bartels_72_01}.  Using the generalized Schur decomposition, we  simultaneously quasi-upper triangularize\footnote{A quasi-upper triangular is  $2 \times 2$ block upper triangular, though generically the blocks on the diagonal are also upper triangular.} $\tilde M_n$ and $\tilde S_n$ via unitary matrices $Q$ and $Z$ that satisfy
	$Q U Z^\top = \tilde M_n$ and $Q T Z^\top = \tilde S_n.$
Thus, writing $Y = X_n^{(2)} Z^\top$  (which is still $\infty \times n$) we have
	$$\L Y U^\top + \N Y  T^\top =  \tilde F_n Q, \qquad {\cal B}_x Y = \vc g_x Q.$$
For simplicity, assume that $U$ and $T$ are upper triangular (the adaption to quasi-upper triangular can be found in~\cite{Bartels_72_01}).  By multiplying the equation by $\vc e_n$, we observe that  the last column of $Y=\bvect[\vc y_1,\dots,\vc y_n]$ satisfies
	$$(U_{nn} \L + T_{nn} \N) \vc y_n  = \tilde F_n Q \vc e_n, \qquad {\cal B}_x \vc y_n = \vc g_x.$$
This equation has the form considered in the previous section and hence is solvable in $\O(n_x^{\rm opt})$ operations with the adaptive QR method (assuming that the sub-problem converges according to the convergence criteria).  The next column satisfies ${\cal B}_x \vc y_{n-1} = \vc g_x Q\vc e_{n-1}$ and 
	\meeq{
		(U_{(n-1)(n-1)} \L + T_{(n-1)(n-1)} \N) \vc y_{n-1}  = \cr
		\tilde F_n Q \vc e_{n-1} - \br[U_{(n-1)n} \L + T_{(n-1)n} \N]\vc y_n.
		}
Thus, we can also calculate $\vc y_{n-1}$.  The procedure continues, calculating each column of $Y$ in turn. Afterwards, we recover $X_n^{(2)} = Y Z$, and then $X_n^{(1)} = \vc g_y - X_n^{(2)} B_n^{(2)\top}$.

\Figurew{pdetime}{\hsize}
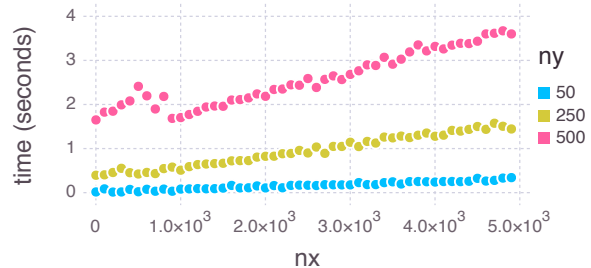
	Computational time to solve Helmholtz equation with homogenous Dirichlet conditions and a forcing term having ${\tt nx} \times {\tt ny}$ ones, ignoring the time for the QZ step.    The complexity is linear in ${\tt nx}$.  It takes less than 4 seconds to solve a PDE with 2.5 million unknowns.

The resulting method has a complexity of $\O(n_y^3 + n_y^2 n_x^{\rm opt})$, where $n_x^{\rm opt}$ is the number of degrees of freedom needed to calculate the worst case column of $Y$.   The $\O(n_y^3)$ term is for the QZ algorithm used to determine the generalized Schur decomposition~\cite{Moler_73_01} and the $\O(n_y^2 n_x^{\rm opt})$ term is for the modified Bartels--Stewart algorithm.   If $n_y\ll n_x$ then this approach has a lower complexity than the fully discretized approach of~\cite{Townsend_14_03}, which achieved $\O(n_y^3 + n_x^3)$ complexity for fixed $n_y$ and $n_x$.    
%When we consider the adaptive solution of PDEs, choosing convergence adaptively  produces considerable savings over the adaptive approach of~\cite{Townsend_14_03}, which progressively doubles $n_x$ until convergence is observed\attodo{I don't understand the $n_x^{opt}$ ``built into'' part of this sentence.}.
  In \figref{pdetime}, we plot the time to solve Helmholtz's equation 
  \meeq{u_{xx} + u_{yy} + 100u = \sum_{k=0}^{{\tt nx}-1} \sum_{j=0}^{{\tt ny}-1} T_k(x) T_j(y), \ccr
  	 u(\pm 1,y)=u(x,\pm 1) = 0,
	 }on $[-1,1]^2$, where $T_k(x) = \cos(k \cos^{-1} x)$ are Chebyshev polynomials. We use the ultraspherical method to represent the partial differential operator as a tensor product of banded operators.  

\begin{remark}
	If the operator is not of splitting rank 2, then in the discrete setting the approach of \cite{Townsend_14_03} is to represent the generalized Sylvester equation via a Kronecker product of the underlying matrices.  This approach extends to the semi-discrete equations as well by interlacing the entries, but results in a substantially higher computational complexity.  
%	though this results in a almost-banded operator with bandwidth $\O(n_y M)$, where $M$ is the largest bandwidth, and hence the complexity reduces to $\O(M^2 n_y^2 n_x^{\rm opt})$, but now we expect $n_x^{\rm opt}$ to be on the order of $n_x n_y$, giving $\O(M^2 n_y^3 n_x)$ operations.  
\end{remark}

%General tensor product equations, after discretizing 
%

\section*{Acknowledgments}
 We acknowledge the support of the Australian Research Council through the Discovery Early Career Research Award (SO).

\end{document}